\documentclass[lineno]{biometrika}
\usepackage[pdftex]{graphicx}
\usepackage{epstopdf}
\usepackage{amsmath,float,url,subfigure,epsfig}
\usepackage{times}
\usepackage{bm}
\usepackage{natbib}
\usepackage[plain,noend]{algorithm2e}

\usepackage[font=small]{caption}

\DeclareMathAlphabet{\mathpzc}{OT1}{pzc}{m}{it}
\usepackage{amssymb,amsfonts}
\usepackage{array,xr}
\usepackage{booktabs,multirow}
\usepackage{graphicx,epsfig,epstopdf}
\usepackage[usenames]{color}
\usepackage{colortbl}

 \makeatletter
\renewcommand{\algocf@captiontext}[2]{#1\algocf@typo. \AlCapFnt{}#2} 
\def\@algocf@capt@plain{top}
\renewcommand{\algocf@makecaption}[2]{%
  \addtolength{\hsize}{\algomargin}%
  \sbox\@tempboxa{\algocf@captiontext{#1}{#2}}%
  \ifdim\wd\@tempboxa >\hsize
    \hskip .5\algomargin%
    \parbox[t]{\hsize}{\algocf@captiontext{#1}{#2}}
  \else%
    \global\@minipagefalse%
    \hbox to\hsize{\box\@tempboxa}
  \fi%
  \addtolength{\hsize}{-\algomargin}%
}
\makeatother
\def\references{\bibliography{1-1-19}}
\newcommand{\be}{\begin{eqnarray}}
\newcommand{\ee}{\end{eqnarray}}
\newcommand{\bea}{\begin{eqnarray*}}
\newcommand{\eea}{\end{eqnarray*}}
\newcommand{\ed}{\end{document}}
\newcommand{\no}{\noindent}
\newcommand{\btab}{\begin{tabular}}
\newcommand{\etab}{\end{tabular}}
\newcommand{\bc}{\begin{center}}
\newcommand{\ec}{\end{center}}

\newcommand{\la}{\label}

\newcommand{\bfi}{\begin{figure}}
\newcommand{\efi}{\end{figure}}
\newcommand{\ben}{\begin{enumerate}}
\newcommand{\een}{\end{enumerate}}
\newcommand{\bdes}{\begin{description}}
\newcommand{\edes}{\end{description}}
\newcommand{\bay}{\begin{array}}
\newcommand{\eay}{\end{array}}

\def\bco{\iffalse}
\def\cov{{\rm cov}}

\def\var{{\rm var}}

\def\ci{\cite}
\def\cp{\citep}

\def\cov{{\rm cov}}

\def\var{{\rm var}}

\def\trace{{\rm trace}}

\def\ci{\cite}
\def\cp{\citep}

\def\eps{\varepsilon}

\def\o{\omega}
\def\O{\Omega}

\def\F{Fr\'{e}chet}

\def\lra{\longrightarrow}

\def\d{\rm d}

\def\V{V_F}
\def\vp{V_p}
\def\hV{\hat{V}_F}

\def\d2{d_2}

\def\hmu{\hat{\mu}_F}
\def\muF{\mu_F}

\def\VF{\hat{V}_F}

\def\s1n{\sum_{i=1}^n}
\def\p1n{\prod_{i=1}^n}
\def\i01{\int_0^1}
\def\1d{{T_\delta}}
\def\1d{{1_\delta}}
\def\1/n{\frac{1}{n}}

\begin{document}

\jname{Biometrika}
\jyear{xxxx}

\markboth{Paromita Dubey \and Hans-Georg M\"uller}{Analysis of Fr\'echet Variance for Metric Space Data}

\title{Fr\'echet Analysis of Variance for Random Objects}
\author{Paromita Dubey}
\affil{Department of Statistics, University of California, Davis,  California 95616, U.S.A. \email{pdubey@ucdavis.edu}}

\author{Hans-Georg M\"uller}
\affil{Department of Statistics, University of California, Davis,  California 95616, U.S.A. \email{hgmueller@ucdavis.edu}
}

\maketitle

\begin{abstract}
		\F \ mean and variance provide a way of obtaining mean and variance for metric space valued random variables and can be used for statistical analysis of  data objects that lie in abstract spaces devoid of algebraic structure and operations. Examples of such data include covariance matrices, graph Laplacians of networks and univariate probability distribution functions. We derive a central limit theorem for \F \ variance under mild regularity conditions,  utilizing empirical process theory,  and also provide a consistent estimator of the asymptotic variance. These results lead to  a  test for comparing $k$ populations of metric space valued data objects in terms of \F \ means and variances. We examine the finite sample performance of this novel inference procedure through simulation studies for several special cases that include probability distributions and graph Laplacians, which leads to a test for comparing  populations of networks. The proposed methodology has good finite sample performance in simulations  for different kinds of random objects.  
		We illustrate the proposed  methods with  data on mortality profiles of various countries  and  resting state Functional Magnetic Resonance Imaging data.
\end{abstract}


\section{Introduction}
\label{sec: intro}

\no With an increasing abundance of complex non-Euclidean data, settings where data objects are assumed to be random variables taking values in a metric space are more frequently encountered. In such settings, where in the general case no manifold or algebraic structure can be assumed, only pairwise distances between the observed data objects are available. The standard problem of k-sample testing is of basic interest in statistics.  For Gaussian data, when comparing means is of primary interest, this corresponds to the classical analysis of variance problem which uses comparisons of sums of squares that measure variation between and within groups. However the k-sample test problem becomes more challenging when data objects lie in general metric spaces.

For general metric space valued random variables, \cite{frec:48} provided a direct generalization of  the mean, which implies a corresponding generalization of variance that may be used to quantify spread of the distribution of metric space valued random variables or random objects around their \F \  mean. The \F \ mean resides in the object space and therefore is not amenable to algebraic operations, which implies that a central limit theorem cannot be directly applied to obtain limit distributions for \F \ means in general metric spaces.  In contrast, the \F \ variance is always a scalar, which makes it more tractable. One of our key results is a central limit theorem for the empirical \F \ variance of data objects in general metric spaces under weak assumptions. While this result is of interest in itself, we proceed to demonstrate how it can be applied to derive a $k$-sample test for comparing \F \ means and variances of population distributions for the case of metric space valued data objects in a spirit similar to the classical analysis of variance problem.
 
In recent years the study of nonparametric tests for the equality of two distributions for Euclidean data has broadened to cover non-Euclidean data, which are increasingly encountered in settings that feature large and complex data. Major approaches have been based on nearest neighbors \cp{henz:88,henz:99,schi:86}, graphs \cp{frie:79,chen:16, rose:05}, energy statistics \cp{szek:04,lyon:13,szek:17} and related work \cp{bari:04}, as well as  kernels \cp{gret:12,sejd:13}. An extension of the energy test for spaces admitting a manifold stratification \cp{patr:15} has been explored recently \citep{guo:17}. 

Empirically, many tests have good power performance for either location type alternatives or scale type alternatives, but usually not for both simultaneously. A major challenge for some of these tests, especially for the case of complex data, is the choice of the required tuning parameters, which often has a major impact on the resulting inference. Some of these challenges associated with existing inference procedures motivate our proposed test, which is simple, easy to compute and mimics the statistic on which classical analysis of variance is based, replacing between and within sums of squares with the corresponding \F \ variances
for separate and combined samples.

\F \ mean based testing and corresponding large sample theory including laws of large numbers and central limit theorems for empirical \F \ means have been explored previously for data objects that lie in a special type of  metric space, such as smooth Riemannian manifolds \cp{bhat:03, bhat:05, bhat:12} and topologically stratified spaces under certain restrictions, like phylogenetic trees \cp{kend:11,bard:13,bhat:13}. Virtually all of these results depend on local linear tangent or similar approximations that are specific to the finite dimensional manifold spaces that are considered in these approaches, and thus rely on local Euclidean approximations. This means 
{existing results on \F \ mean and variance testing} do not apply for random objects in more general {infinite dimensional} metric spaces, such as the space of probability density functions. The central limit theorem for \F \ means was recently applied to the space of graph Laplacians \cp{gine:17}, which are of interest to obtain inference for networks. This required choosing a high dimension for the approximating space, thus leading to problems with small sample high dimensional data and the ensuing complications for inference.

Our goal in this paper is to develop a simple and straightforward extension of analysis of variance to the case of metric space valued random variables. Our starting point is a bounded metric space $(\O,d)$. We show that consistency of the sample \F$\,$ mean can be obtained  by using results of \cite{mull:19:3} concerning \F \ regression estimators. We derive a central limit theorem for \F \ variance and provide a consistent estimator of its asymptotic variance. Our method is applicable to a wide class of objects including correlation matrices,  probability distributions, manifolds and also the space of graph Laplacians. Making use of this new central limit theorem, we derive a  $k$-sample test for random objects and study the asymptotic distribution of the test statistic under the null hypothesis of equality of the population distributions, as well as its power function.

It is customary to test for heteroscedasticity of the population groups prior to applying an $F$-test in classical analysis of variance, in order  to evaluate the assumption of equal variances across the populations that are compared. One popular test for this purpose \citep{leve:60} is of the form of the usual analysis of variance $F$-test, but applied to pseudo-observations which could in principle be any monotonic function of the absolute deviations of the observations from their group ‘centers’. In our approach, Levene's  test and the classical analysis of variance for testing inequality of population means are combined to derive a test for data objects in general metric spaces. This makes it possible to aim at both location and scale type alternatives, instead of only location alternatives as in classical analysis of variance. 

In addition to the asymptotic distribution of the test statistic and the consistency of the proposed test, we also derive results that justify a bootstrap implementation, which is shown to be useful for smaller sample sizes in simulations. We demonstrate the implementation of the proposed test in simulations and in applications that include the comparison of networks and of  density functions in demography and of mentally normal and Alzheimer's subjects based on fMRI brain imaging data.

\bco

 Our test statistic is a sum of two components. One component is proportional to the squared difference of the pooled sample \F \ variance and the weighted
average of the groupwise \F \ variances, with weights proportional to the sample sizes of the groups which for the special case of Euclidean data, is proportional to the squared F-ratio as in the usual ANOVA. A key auxiliary result is that this statistic converges to zero at rate $o_P(1/n)$ under the null hypothesis of equality of \F \ means of population distributions. The other component of our test statistic accounts for differences in the \F \ variances of the population groups and under the Euclidean setting simplifies to a generalization of the Levene’s test applied to squared absolute deviations of the observations from their group \F \ means. It turns out that when the assumptions of the central limit theorem hold and under the null hypothesis of equality of \F \ means and variances of the populations, the asymptotic distribution of our test statistic is $\chi^2_{k-1}$.

\fi

\section{Preliminaries}
\label{sec: prelim}
\no The \F \ mean is a generalization of centroids to metric spaces and for the special case of  Euclidean data it includes the arithmetic mean, median and geometric mean under different choices of distance functions. The \F \ variance is the corresponding  generalized measure of dispersion around the \F \ mean. More formally, in all of the following, $(\O, d, P)$ is a bounded metric space with metric $d$ and probability measure $P$. Random objects  are  random variables $Y$ that take values in $\O$.  The population \F \ mean $\mu_F$ of $Y$ and the sample \F \ mean $\hmu$ for  a random sample  $Y_1, Y_2,\dots, Y_n$ of independent and identically distributed random variables with the same distribution  as $Y$ are given by 
\begin{equation*}
\label{eq: fre_mean_pop}
\mu_F=\operatornamewithlimits{argmin}_{\omega \in \Omega}E\{d^2(\omega, Y)\}, \quad \hmu=\operatornamewithlimits{argmin}_{\omega \in \Omega} \frac{1}{n}\s1n d^2(\omega,Y_i).
\end{equation*}

The sample \F \ mean is an $M$-estimator as it is obtained by minimizing a sum of functions of the data objects. The \F \ variance quantifies the spread of the random variable $Y$ around its \F \ mean $\mu_F$. The population \F \ variance $\V$ and its sample version $\VF$ are 
\begin{equation*}
\label{eq: fre_var_pop}
V_F= E\{d^2(\muF,Y)\}, \quad \hV=\frac{1}{n}\sum_{i=1}^{n}d^2(\hmu,Y_i).
\end{equation*}
Note that $\mu_F, \hmu \in \O$, while $V_F, \hV \in \mathbb{R}$. The asymptotic consistency of the sample \F \ mean $\hmu$ follows from  results in \cite{mull:19:3}, under the following assumption:
\begin{assumption}
	The objects $\hat{\mu}_F$ and $\mu_F$ exist and are unique, and  for any $\epsilon > 0$, \ $\inf_{d(\omega,\mu_F)>\epsilon} E\{d^2(\omega,Y)\}> E\{d^2(\mu_F,Y)\}$.
\end{assumption}
Assumption 1 is instrumental to 
establish the weak convergence of the empirical process $H_n(\o)=\1/n \s1n d^2(\o,Y_i)$ to the population process $H(\o)= E \{d^2(\o,Y)\}$, which implies the consistency of $\hat{\mu}_F$, 
\begin{equation}
\label{eq: fre_mean_consistency}
d(\hat{\mu}_F,\mu_F) = o_P(1).
\end{equation}
Observing that 
\begin{align}
\label{eq: fre_var_consistency}
|\VF-V_F| \leq \left| \frac{1}{n} \sum_{i=1}^n \left \lbrace d^2(Y_i,\hat{\mu}_F)-d^2(Y_i,\mu_F) \right \rbrace \right| +  \left| \frac{1}{n} \sum_{i=1}^n d^2(Y_i,{\mu}_F)-V_F \right|,
\end{align}
consistency of $\hmu$ is  seen to imply the consistency of $\VF$. This is because the first term in \eqref{eq: fre_var_consistency} is upper bounded by $2 \text{diam}(\O) d(\hat{\mu}_F,\mu_F) $ where $\text{diam}(\Omega)=\sup \{d(\omega_1,\o_2): \o_1, \o_2 \in \O\}$ is finite since $\O$ is bounded, and the second term in \eqref{eq: fre_var_consistency} is $o_P(1)$ by the weak law of large numbers.

For a central limit theorem to hold for the empirical \F \ variance we need an assumption  on the complexity of the metric space $\O$, which can be quantified by a bound on the entropy integral for metric $\delta$-balls $B_\delta(\o)$ 
of  $\O$ \cp{well:96}, given by 
\begin{equation*}
J(\delta,\o)=\int_{0}^{1} \left[1+\log N \{\epsilon\delta/2,B_\delta(\o),d\} \right]^{1/2} \ d\epsilon,
\end{equation*} 
where $B_\delta(\o)$ is the $\delta$-ball in the metric $d$, centered at $\o$ and $N\{\epsilon\delta/2,B_\delta(\o),d\}$ is the covering number for $B_\delta(\o)$ using open balls of radius $\epsilon\delta/2$. Specifically, to obtain the desired central limit theorem we assume that 
\begin{assumption}
 For any $\o \in \O$,	$\delta J(\delta,\o) \rightarrow 0$ as $\delta \rightarrow 0$.
\end{assumption}
For our results on the power of the proposed test in Section 4, we need an additional assumption on the entropy integral of the whole space $\O$.
\begin{assumption}
	The entropy integral of  $\O$ is finite,  $\int_{0}^{1} \left \{1+\log{N(\epsilon,\O,d)}\right \}^{1/2}d\epsilon < \infty$.
\end{assumption}	
Examples of random objects include univariate probability distributions with supports that are contained in a common compact set  in $\mathbb{R}$ with finite second moments, equipped with the Wasserstein metric $d_W$, and also the  spaces of correlation matrices and graph Laplacians of fixed dimensions of weighted networks with bounded weights, equipped with the Frobenius metric $d_F$.  {In all these examples, the metric space $(\O,d)$ is bounded.}  For two univariate distributions $F$ and $G$ with finite variances, the $L^2$-Wasserstein distance, also known as earth movers distance, which is connected to optimal transport  \cp{vill:03},  is given by $d^2_W(F,G)= \int_{0}^{1} \{F^{-1}(t)-G^{-1}(t)\}^2dt$, where $F^{-1}$ and $G^{-1}$ are the quantile functions corresponding to $F$ and $G$, respectively. For two matrices $A$ and $B$ of the  same dimension, we consider the  Frobenius metric given by $d^2_F(A,B)= \trace\{(A-B)'(A-B)\}$. According to Propositions 1 and 2 in \cite{mull:19:3}, the space $(\O,d_W)$ satisfies Assumptions 1-3 when the set $\O$  consists  of univariate probability distributions {with compact support} in $\mathbb{R}$ and having finite second moments and $d_W$ is the $L^2$-Wasserstein metric, as does the space $(\O,d_F)$ when the set $\O$  consists of correlation matrices of a fixed dimension and $d_F$ is the Frobenius metric.

To characterize the space of graph Laplacians, we denote a weighted undirected graph by $G=(V,E)$, where $V$ is the set of its vertices and $E$ the set of its edges. Given an  adjacency matrix $W$, where $0 \leq w_{ij}=w_{ji} \leq 1$ and equality with zero holds if and only if $\{i,j\} \notin E$, the graph Laplacian is defined as  $L=D-W$, where $D$ is the diagonal matrix of the degrees of the vertices, i.e. $d_{jj}=\sum_{i}w_{ij}$. Under the assumption that the graphs are simple, i.e. there are no self loops or multi edges, there is a one to one correspondence between the space of graphs and the graph Laplacians. Hence the space of graph Laplacians can be used to characterize the space of networks \cp{gine:17}. It is easy to see that the space $(\O,d_F)$ satisfies Assumptions 1-3 when the set $\O$  consists of graph Laplacians of connected, undirected and simple graphs of a fixed dimension by a minor extension of the arguments provided for the space of correlation matrices of a fixed dimension in \cite{mull:19:3}.

Assumptions 2 and 3 help control the complexity of $(\O,d)$ 
through the covering numbers $N(\eps,\O,d)$ and $N\{\eps\delta,B_\delta(\o),d\}$ by restricting their rates of increase as  $\epsilon > 0$ decreases. An example of a slightly stronger but sometimes  more interpretable condition, as compared to assumptions 2 and 3, is the upper bound on  the metric entropy of the space $(\O,d)$ given by 
	\begin{equation*}
	\log N(\eps,\O,d) \leq K/\eps^\alpha,
	\end{equation*}
for some constants $K>0$ and $\alpha < 2$.  This stronger condition is satisfied by a wide class of metric spaces, including the space of distribution functions on $\mathbb{R}$ and  $\mathbb{R}^2$ with the $L^2$ metric, a consequence of Theorem 2.6.9 in \citep{well:96}, the space of  monotone functions from $\mathbb{R}$ to a compact subset of $\mathbb{R}$, the space of all Lipschitz functions of degree $\gamma \leq 1$ on the unit interval $[0,1]$ with the $L^2$ metric, as well as the class of convex functions on a compact, convex subset of $\mathbb{R}^d$ under certain restrictions \cp{well:96}.

\section{Central Limit Theorem for \F \ Variance}
\label{sec: theory}
The following Proposition  lays the foundations for the Central Limit Theorem for the empirical \F \ variance $\VF$.
\begin{proposition}
	\label{lma: main}
	Suppose Assumptions 1-3 hold. Then 
	\begin{equation*}
	\label{eq: main}
	\frac{1}{n}\sum_{i=1}^{n}\{d^2(\hat{\mu}_F,Y_i)-d^2(\mu_F,Y_i)\}  = o_P(n^{-1/2}).
	\end{equation*}
\end{proposition}
All proofs are in the Supplementary Material. Proposition  \ref{lma: main} makes it possible to  deal with the sum of dependent random variables $\s1n d^2(\hmu,Y_i)$ by replacing it with the sum of i.i.d. random variables $\s1n d^2(\muF,Y_i)$, which is a crucial step in the derivation of the central limit theorem for $\VF$. 
Since $\O$ is bounded, {the quantity $\var\{d^2(\mu_F,Y)\}$ appearing in Theorem \ref{thm: clt}, which is the asymptotic variance of the \F \ variance}, is always finite. Possible extensions to general M-estimators are discussed in the Supplementary Material.  The  central limit theorem for \F \ variance is as follows. 
\begin{theorem}
	\label{thm: clt}
	Under the assumptions of Proposition  \ref{lma: main},
	\begin{equation*}
	\label{eq: clt}
	{n}^{1/2}(\VF-V_F)\lra N(0,\sigma_F^2) \quad \text{in distribution}, 
	\end{equation*}
	where $\sigma_F^2=\var\{d^2(\mu_F,Y)\}$.
\end{theorem}
An intuitive sample based estimator for $\sigma_F^2$ is  
\begin{equation}
\label{eq: var_est}
\hat{\sigma}_F^2=\frac{1}{n}\sum_{i=1}^{n} d^4(\hat{\mu}_F,Y_i)-\left \{ \frac{1}{n}\sum_{i=1}^{n}d^2(\hat{\mu}_F,Y_i)\right \} ^2  , 
\end{equation}
and $n^{1/2}\left(\widehat{\sigma}_F^2 -\sigma_F^2 \right)$ has an asymptotic normal distribution, as follows. 
\begin{proposition}
	\label{lma: consistency}
	Under the assumptions of Proposition  \ref{lma: main},
	\begin{equation*}
	\label{eq: consistency}
	n^{1/2} \left( \widehat{\sigma}_F^2 -\sigma_F^2 \right) \lra  N(0,A) \quad \text{in distribution}, 
	\end{equation*}
	where $A=a'Ba$ with $a'= \left[1,
	-2E\left\{d^2(\mu_F,Y)\right\}\right]$ and $$B=\begin{bmatrix}
	\var\left\{d^4(\mu_F,Y)\right\} & \cov \left \{d^4(\mu_F,Y),d^2(\mu_F,Y)\right \} \\
	\cov \left \{d^4(\mu_F,Y),d^2(\mu_F,Y)\right \}& \var\left \{d^2(\mu_F,Y)\right \}
	\end{bmatrix}.$$
\end{proposition}

The matrix $A$ is well-defined due to the boundedness of $\O$. Combining Theorem \ref{thm: clt}, Proposition \ref{lma: consistency} and Slutsky's Theorem leads to 
\begin{equation}
\label{eq: ci_1}
n^{1/2}(\VF-V_F) / \hat{\sigma}_F \lra N(0,1) \quad \text{in distribution}. 
\end{equation} 
A simple application of the delta method gives the asymptotic distribution of the \F \ standard deviation, defined as  the square root of the  \F \ variance,			\begin{equation*} \la{std}
n^{1/2}({\VF}^{1/2}-{V_F}^{1/2})\lra N(0,\sigma_F^2/4 V_F) \quad \text{in distribution}, 
\end{equation*}
and since both $\hat{\sigma}_F$ and $\VF$ are consistent estimators,
\begin{equation}
\label{eq: ci_2}
2 n^{1/2}  {\VF}^{1/2} ({\VF}^{1/2}-{V_F}^{1/2})/\hat{\sigma}_F\lra N(0,1) \quad \text{in distribution}. 
\end{equation}

One can use \eqref{eq: ci_1} and \eqref{eq: ci_2} to construct asymptotic confidence intervals for \F \ variance and standard deviation, which depend on the quality of the large sample approximations. The bootstrap provides an alternative that often has better finite sample properties under weak assumptions 
\cp{bick:81,bera:03}. Under fairly general assumptions, resampling methods like bootstrapping and permutation tests work whenever a central limit theorem holds \cp{jans:03}.  A basic criterion for bootstrap confidence sets to have correct coverage probability asymptotically is convergence of the bootstrap distribution of the root, in our case
$n^{1/2}({\VF}-{V_F})/{\hat{\sigma}_F}$.Then Monte Carlo approximations of the bootstrap distribution of the root provide  approximate quantiles for the construction of confidence sets. Further details are provided in the Supplementary Material. We conclude  that the nonparametric bootstrap is a viable option for the construction of confidence intervals for Fr\'echet variance. 

\section{Comparing Populations of Random Objects}
\label{sec: k-sample test}
 \no  Assume we have  a  sample of $\O$-valued random data objects $Y_1,Y_2,\dots,Y_n$ 
that belong to $k$ different groups $G_1,G_2,\ldots,G_k$, each of size $n_j \, (j=1,\ldots,k),$ such that
$\sum_{j=1}^{k} n_j=n$. We wish to test the null hypothesis that the \F \ means and variances of the population distributions of the $k$ groups are identical versus the alternative that at least one of the groups has a different population distribution compared to the others in terms of either its \F \ mean or its \F \ variance. Consider  the sample \F \ means 
\begin{equation*}
\hat{\mu}_j=\operatornamewithlimits{argmin}_{\omega \in \Omega} \frac{1}{n_j}\sum_{i \in G_j}d^2(\omega,Y_i) ,
\end{equation*} which are random objects computed just from the data falling into group $j$. The  corresponding real-valued sample \F \ variances are 
\begin{equation*}
\hat{V}_j=\frac{1}{n_j}\sum_{i \in G_j} d^2(\hat{\mu}_j,Y_i), 
\end{equation*} with associated  variance estimates (\ref{eq: var_est}) given by 
\begin{equation*}
\hat{\sigma}^2_j=\frac{1}{n_j}\sum_{i \in G_j} d^4(\hat{\mu}_j,Y_i) - \left \{\frac{1}{n_j}\sum_{i \in G_j} d^2(\hat{\mu}_j,Y_i) \right\}^2 \quad (j=1,\ldots,k).
\end{equation*}Consider also the  pooled sample \F \ mean $\hat{\mu}_p$ and its corresponding pooled sample \F\ variance  $\hat{V}_p$, 
\be \la{pe} \hat{\mu}_p=\operatornamewithlimits{argmin}_{\omega \in \Omega} \frac{1}{n}\sum_{j=1}^{k}\sum_{i \in G_j} d^2(\omega,Y_i), \quad \hat{V}_p= \frac{1}{n}\sum_{j=1}^{k}\sum_{i \in G_j} d^2(\hat{\mu}_p,Y_i),  \ee
and 
define  $\lambda_{j,n}=n_j/n \,\,  (j=1,\dots,k),$ where $\sum_{j=1}^{k} \lambda_{j,n}=1$.

We  will base our inference procedures on the auxiliary statistics
\begin{equation}
\label{eq: factor}
F_n= \hat{V}_p-\sum_{j=1}^{k}\lambda_{j,n} \hat{V}_j
\end{equation}
and 
\begin{equation}
\label{eq: var_stat}
U_n=  \sum_{j < l} \frac{\lambda_{j,n}\lambda_{l,n}}{\hat{\sigma}_j^2 \hat{\sigma}_l^2} (\hat{V}_j-\hat{V}_l)^2.
\end{equation}
Here $F_n$ is almost surely non-negative and is equal to the  numerator of the $F$-ratio in classical Euclidean analysis of variance,  where it corresponds to the weighted variance of the group means, with weights proportional to the group sizes.  
It can be regarded as a generalization of the F-ratio in classical ANOVA to the more general setting of metric space valued data. Analogous to the classical scenario, $F_n$ is expected to be small under the null hypothesis of equality of \F \ means of the population distributions, which is indeed the case as demonstrated by the following Proposition.
\begin{proposition}
	\label{lma: f-ratio}
	Suppose $\hat{\mu}_p$ and  $\hat{\mu}_j$ exist and are unique almost surely for all $j=1,\ldots,k$. Let $0  < \lambda_{j,n} < 1 \,\, (j=1,\dots,k)$ and $\lambda_{j,n} \rightarrow \lambda_j \ \text{as} \ n \rightarrow \infty$, where $\lambda_j$ is such that $0 < \lambda_j < 1 \quad (j=1,\dots,k),$ with $\sum_{k=1}^k \lambda_j =1.$ Then under the null hypothesis of equality of \F \ means of the population distributions and under Assumptions 1 and 2 for each of the groups, as $n \rightarrow \infty$, 
	\begin{equation}
	\label{eq: f-ratio}
	n^{1/2} F_n=o_P(1).
	\end{equation}
\end{proposition}
Inference in classical ANOVA requires Gaussianity and 
equality of the population variances and hence targets  only differences in the group means 
to capture differences in the population distributions. 
Aiming at detecting a broader class of alternatives, we  employ the statistics 
$U_n$ in \eqref{eq: var_stat} to target differences among the population variances, where in 
the Euclidean case,  $U_n$ turns out to be a slightly modified version of the traditional Levene's test, substituting squared distances of the observations from their group \F \ means instead of just distances. The following result provides the asymptotic distribution of $U_n$ under the null hypothesis of equal population \F \ variances.
\begin{proposition}
	\label{lma: u-stat}
	Under the assumptions of Proposition  \ref{lma: f-ratio}, we have under the null hypothesis of equal population \F \ variances, as $n \rightarrow \infty$, 	    	
	\begin{equation}
	\label{eq: u-stat}
	\frac{nU_n}{\sum_{j=1}^{k}\frac{\lambda_{j,n}}{ \hat{\sigma}_j^2}} \lra \chi^2_{(k-1)} \quad \text{in distribution}. 
	\end{equation}
\end{proposition}
The proposed test statistic $T_n$ is then  
\begin{equation}
\label{eq:test_stat}
T_n= \frac{n U_n }{\sum_{j=1}^{k}\frac{\lambda_{j,n}}{ \hat{\sigma}_j^2}} + \frac{nF_n^2}{\sum_{j=1}^{k}\lambda_{j,n}^2\hat{\sigma}_j^2}.
\end{equation}
When constructing $T_n$, we scale $F_n$ by an estimate of the standard deviation of $\sum_{j=1}^{k}\lambda_{j,n}\hat{V}_j$, which is $(\sum_{j=1}^{k}\lambda_{j,n}^2\hat{\sigma}_j^2)^{1/2}$,  so that $F_n$ is suitably scaled with respect to the variability of $\sum_{j=1}^{k}\lambda_{j,n}\hat{V}_j$,  and then square it so that both terms in $T_n$ are of the same order in $n$. We place equal weights on both terms in $T_n$ as we assume that there is no prior information on whether differences in population distributions arise due to inequality of their \F \ means or \F \ variances. If such prior information is available, one can consider a modified version that corresponds to a convex average of the two terms in (\ref{eq:test_stat}). 

Under  the null hypothesis of equality of \F \ means, consistency of $\hat{\sigma}_j^2 \,\, (j=1,\ldots,k)$ and Proposition \ref{lma: f-ratio} imply that  
$nF_n^2/\sum_{j=1}^{k}\lambda_{j,n}^2\hat{\sigma}_j^2=o_P(1)$. In combination with Proposition \ref{lma: u-stat}, this  leads to the following result.
\begin{theorem}
	\label{cor: test_stat}
	Under the null hypothesis of equal population \F \ means and variances and the assumptions of Proposition  \ref{lma: f-ratio},
	\begin{equation}
	\label{eq: null-dist}
	T_n \lra \chi^2_{(k-1)} \quad \text{in distribution}. 
	\end{equation}
\end{theorem}

For a level $\alpha$ test, we accordingly  reject the null hypothesis of simultaneous equality of \F \ means and variances if the test statistic $T_n$ turns out to be bigger than $\chi^2_{k-1,\alpha}$, which is the $(1-\alpha)^{th}$ quantile of the $\chi^2_{(k-1)}$ distribution, i.e. the rejection region that defines the  test is 
\be  R_{n,\alpha}=\{T_n > \chi^2_{k-1,\alpha}\}. \la{rej} \ee	     

To study  the consistency of the proposed test  (\ref{rej}), we consider contiguous alternatives that  capture departures from the null hypothesis of equality of  \F \ means and variances. Consider the following population quantities, 
\begin{equation}
\label{eq: pooled}
\mu_{p}= \operatornamewithlimits{argmin}_{\omega \in \Omega} \sum_{j=1}^{k}
\lambda_{j} E_j\{d^2(\omega,Y_j)\}, \hspace{0.2 in}
\vp=  \sum_{j=1}^{k}
\lambda_{j} E_j\{d^2(\mu_{p},Y_j)\}
\end{equation} and
\begin{equation}
F=\vp-\sum_{j=1}^{k}\lambda_{j} V_j , \hspace{0.2 in} U= \ \sum_{j < l} \frac{\lambda_{j}\lambda_{l}}{\sigma_j^2 \sigma_l^2} (V_j-V_l)^2,
\end{equation} where $E_j(\cdot)$ denotes expectation under the probability distribution  for the $j$th population, $\lambda_{j,n}$ is as defined in Proposition \ref{lma: f-ratio}, and the $Y_j$ are random objects distributed according to the  $j$th  population distribution. In  the Euclidean setting, $\mu_p$ and $V_p$ are analogous to the pooled population \F \ mean and pooled population \F \ variance, respectively,  and in the general case can be interpreted as generalizations of these quantities. 

Proposition \ref{prop: pooled} below states that under mild assumptions on the existence and uniqueness of the pooled and the groupwise \F \ means, the statistics $\hat{\mu}_p$ and $F_n$ are  consistent estimators of the population quantities $\mu_p$ and $F$. By our assumptions, $F$ is zero only under the equality of the population \F \ means and positive otherwise. The population quantity $U$ is proportional to the weighted average of the pairwise differences between the groupwise \F \ variances, which is nonnegative and is zero only if the population groupwise \F \ variances are all equal. 
\begin{proposition}
	\label{prop: pooled}
	Suppose $\hat{\mu}_p$, $\hat{\mu}_j$, $\mu_p$ and $\mu_j$ exist and are unique, the sample based estimators almost surely,  for all $j=1,\ldots,k$. Assume that for any $\epsilon > 0,$ $\inf_{d(\omega,\mu_p)>\epsilon}  \sum_{j=1}^{k}
	\lambda_j E_j\{d^2(\omega,Y_j)\} > \sum_{j=1}^{k}
	\lambda_j E_j\{d^2(\mu_{p},Y_j)\}$ and also that  $\inf_{d(\omega,\mu_j)>\epsilon} E_j\{d^2(\omega,Y_j)\} > E_j\{d^2(\mu_j,Y_j)\}$ for all $j=1,\ldots,k$. Let $0  < \lambda_{j,n} < 1 \,\,  (j=1,\dots,k)$ and $\lambda_{j,n} \rightarrow \lambda_j  \ \text{as} \ n \rightarrow \infty$, where $\lambda_j$ is such that $0 < \lambda_j < 1 \,\, (j=1,\dots,k)$, with $\sum_{k=1}^k \lambda_j =1,$ as defined in Proposition \ref{lma: f-ratio}. Then, as $n \rightarrow \infty$,
	\begin{equation}
	d(\hat{\mu}_p , \mu_p)=o_P(1) \quad \ \text{and} \quad |F_n-F| =o_P(1).
	\end{equation}
	The population quantity F is nonnegative and is zero if and only if the population \F \ means $\mu_j$ are all equal. 
\end{proposition}

To study the power performance of the proposed  test (\ref{rej}), we consider sequences of alternatives $H_n$ where $H_n=\{(U,F): F \geq a_n \ \text{or} \  U \geq b_n\}$ for nonnegative sequences $a_n$ or $b_n$ with either $a_n$ or $b_n$ strictly greater than 0. The case where either $a_n \rightarrow 0$ or $b_n \rightarrow 0$, as $n \rightarrow \infty$ reflects contiguous alternatives. Of interest is the asymptotic behavior of the power function $\beta_{H_{n}}$, where
\be
\beta_{H_{n}}= \inf_{(U,F) \in H_n} P\left( R_{n,\alpha}\right).   \la{pow}     
\ee	
Here the rejection region $R_{n,\alpha}$ is as in (\ref{rej}). The 
following result   provides sufficient conditions for the consistency of the  proposed test for this family of contiguous alternatives, which means that the asymptotic power is 1 for any such alternative and  any choice of the level $\alpha$.  

\begin{theorem}
	\label{prop: power}
	Under the assumptions of Proposition \ref{prop: pooled} and Assumption 3, for sequences of  contiguous alternatives  $H_{n}$ for which either $F \geq a_n$ or $U \geq b_n$, where  $a_n \rightarrow 0$ and $b_n \rightarrow 0$ as $n \rightarrow \infty$,   the power function (\ref{pow}) satisfies for all $\alpha >0$\\   
	(A) If $n^{1/2} a_n \rightarrow \infty$, then 
	$\beta_{H_{n}} \rightarrow 1$. \\
	(B)  If  $nb_n \rightarrow \infty$ , then $\beta_{H_{n}} \rightarrow 1$.
\end{theorem}

Asymptotic tests may not work very well for common situations where the group sample sizes are small. By the asymptotic justification provided at the end of Section \ref{sec: theory},  
resampling methods such as  bootstrap and permutation tests using the proposed test statistic $T_n$ can then  be used to obtain more accurate level $\alpha$ tests. 
\section{Simulation Studies}
\label{sec: simu}
	In order to gauge the performance of the proposed test (\ref{rej}), we performed  simulation experiments under various settings. The random objects  we consider  include samples of univariate probability distributions equipped with the $L^2$-Wasserstein metric, samples of graph Laplacians of scale free networks from the Barab\'{a}si-Albert model \cp{bara:99} with the Frobenius metric and samples of multivariate data with the usual Euclidean metric. In each case we considered two groups of equal size $n_1=n_2=100$ and constructed the empirical power functions of the proposed test against departures from the null hypothesis of equality of  population \F \ means and variances of the two groups at level 0.05. The empirical power was computed by finding the proportion of rejections for 1000 Monte Carlo runs. 
	For comparing the performance of the proposed test against existing tests we used the bootstrap version of all the tests, where  the critical values were  obtained from the bootstrap distribution of the test statistics in 1000 Monte Carlo simulations. 
	We also investigated the finite sample power of the proposed asymptotic test for increasing sample sizes by comparing power functions for group sizes $n_1=n_2=100, 250, 450$.

    In the simulations we explored not only location differences but also differences in shape and scale of the population distributions. We compared the proposed  test (\ref{rej})  with the graph based test \citep{chen:16} and the energy test based on pairwise distances \citep{szek:04} 
    For the graph based test, we constructed the similarity graph of the pooled observations of the two groups by constructing  a $5-$MST graph from the pooled pairwise distance matrix, following the suggestion in \ci{chen:16}, where MST stands for minimal spanning tree. Here a  $k-$MST is the union of the $1^{st},\dots, k^{th}$ MST, where the $k^{th}$ MST is a spanning tree connecting all observations that minimizes the sum of distances across edges subject to the constraint that this spanning tree does not contain any edge in the $1^{st},\dots , (k-1)^{th}$  MST. For computing the statistic of the energy test of \cite{szek:04}, we used the pairwise distance matrix obtained from the specified metric in the space of random objects. 

\begin{figure}
	\centering
	   \includegraphics[width=0.85\linewidth]{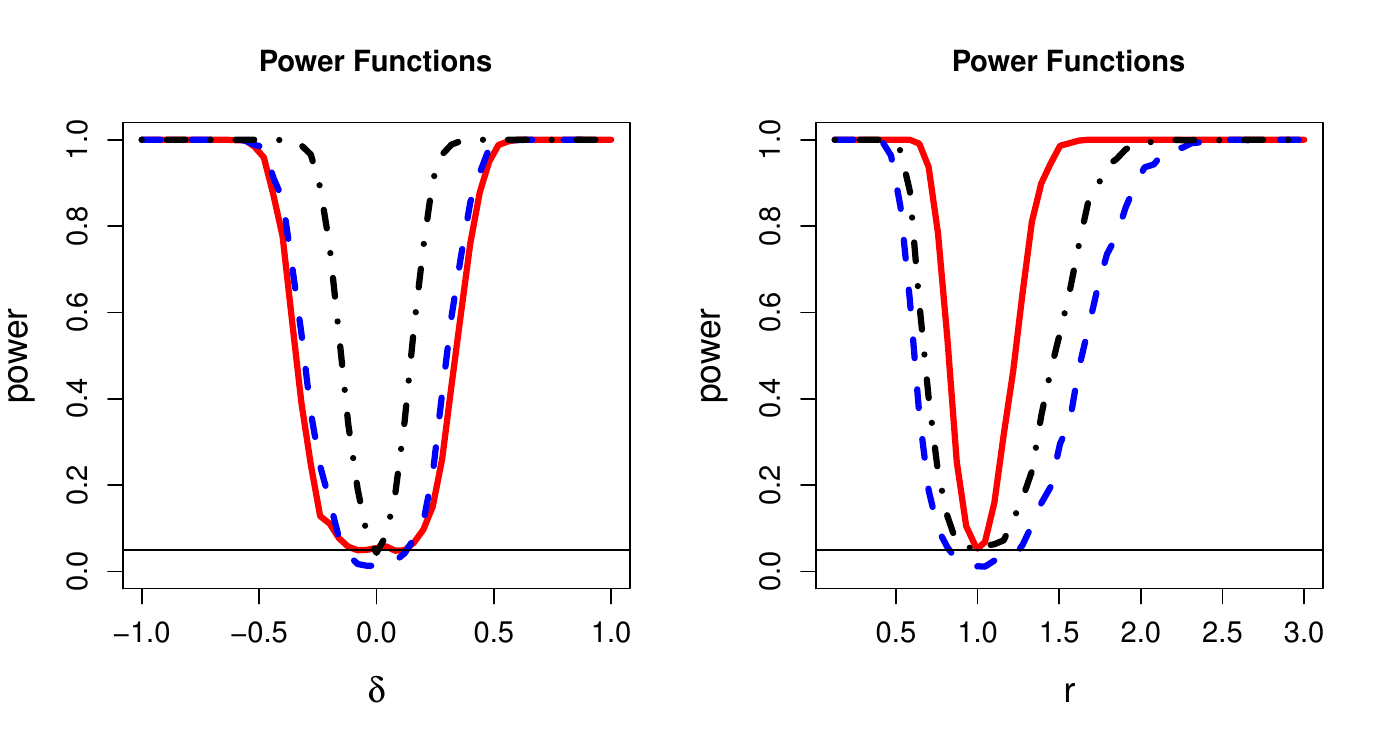}
		\caption{Empirical power as function of  $\delta$ for  $N(\mu,1)$ probability distributions with $\mu$ from $N(0,0.5)$ for group $G_1$ and $N(\delta, 0.5)$  for group $G_2$, {truncated to lie in $[-10,10]$ for both groups} (left) and empirical power as function of  $r$ for  $N(\mu,1)$ probability distributions with $\mu$ from $N(0,0.2)$ for $G_1$ and $N(0,0.2r)$  for $G_2$, {also truncated to lie in $[-10,10]$ for both groups} (right).  The solid red curve corresponds to the bootstrapped version of the proposed test (\ref{rej}),  the dashed blue curve to the graph based test \cp{chen:16} and the dot-dashed black curve to the energy test \cp{szek:04}.  
			The level of the tests is $\alpha= 0.05$ and is indicated by the line parallel to the $x$ axis. Sample sizes of the groups are fixed at $n_1=n_2=100$.}
		\label{fig:fig_1}
    \end{figure}

   The first type of random objects we study are random samples of univariate probability distributions. Each datum is a $N(\mu,1)$ distribution where $\mu$ is random. As distance  between two probability distributions we choose the $L^2$-Wasserstein metric. In the first scenario, for group $G_1$, we generate  $\mu$ {to be distributed as truncated normal distribution $N(0,0.5)$, constrained to lie in $[-10,10]$ and for group $G_2$  as $N(\delta,0.5)$, truncated within $[-10,10]$ and compute the empirical power function of the tests for $-1 \leq \delta \leq 1$}. In the second scenario $\mu$ is drawn randomly from $N(0,0.2)$ for group $G_1$ and from $N(0,0.2r)$ for group $G_2$, {truncated within $[-10,10]$ in both cases,} and empirical power is evaluated for $0.125 \leq r \leq 3$. The first  scenario emphasizes  location differences between the populations and the second emphasizes scale differences. The results are presented in Figure \ref{fig:fig_1}. We find that in  the first scenario of mean differences, the proposed test and the graph based test perform similarly. Both are outperformed by the energy test. In the second scenario of scale differences the proposed test outperforms all other tests.  
  
\begin{figure}
	\centering
		\includegraphics[width=0.85\linewidth]{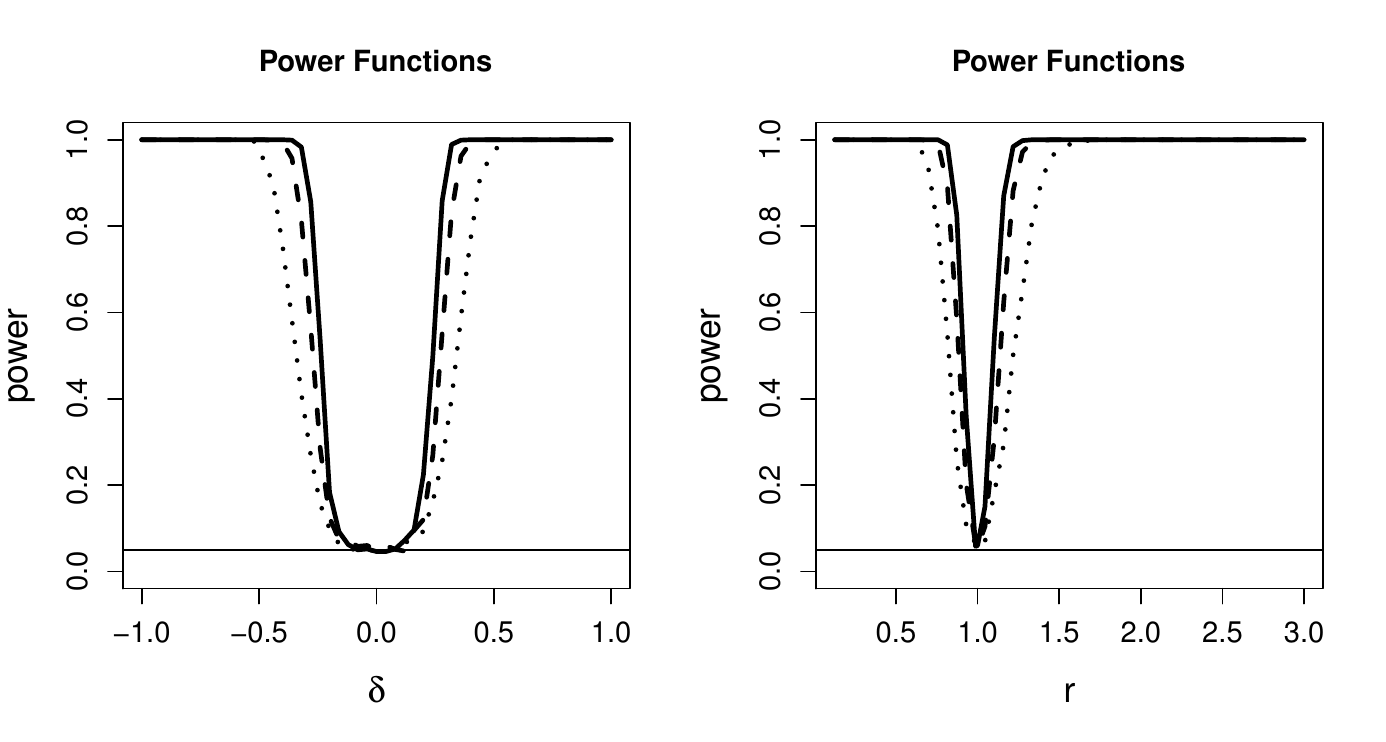}
		\caption{Empirical power in the same setting as in Figure 1 for the proposed test (\ref{rej}) for different sample sizes. 
		The tests are at level $\alpha=0.05$, indicated by the line parallel to the $x$-axis. The  solid curve corresponds 
		to sample sizes	$n_1=n_2=450$, the dashed curve to  $n_1=n_2=250$ and the dotted curve to  $n_1=n_2=100$.}
		\label{fig: fig_2}
\end{figure}

  Theorem \ref{prop: power} provides the theoretical underpinning for  the power analysis of the proposed test, where power against mean differences is controlled by $F$ and against variance differences by $U$. In the simple Euclidean case, the $\mu_j$ are the regular group means. $F$ is then proportional to $\sum_{j=1}^{k} \lambda_i \lambda_j (\mu_i-\mu_j)^2$ while $U$ is proportional to $\sum_{j=1}^{k} \lambda_i \lambda_j (V_i-V_j)^2$, assuming $\sigma_j$ is the same for all the groups. Theorem \ref{prop: power} implies that close to the null, the test gains power when either $F \approx {a_n}/{n^{1/2}}$ or when $U \approx {b_n}/{n}$, which explains why the proposed test performs slightly better in detecting variance differences than mean differences in the proximity of the null hypothesis.
Figure \ref{fig: fig_2} indicates that  the proposed test is consistent for large sample sizes in both  scenarios. {Additional simulations with unequal group sizes in the above setting and comparison between the asymptotic and bootstrap versions of our proposed test are included in the Supplementary Material.}

\begin{figure}
	\centering
		\includegraphics[width=0.85\linewidth]{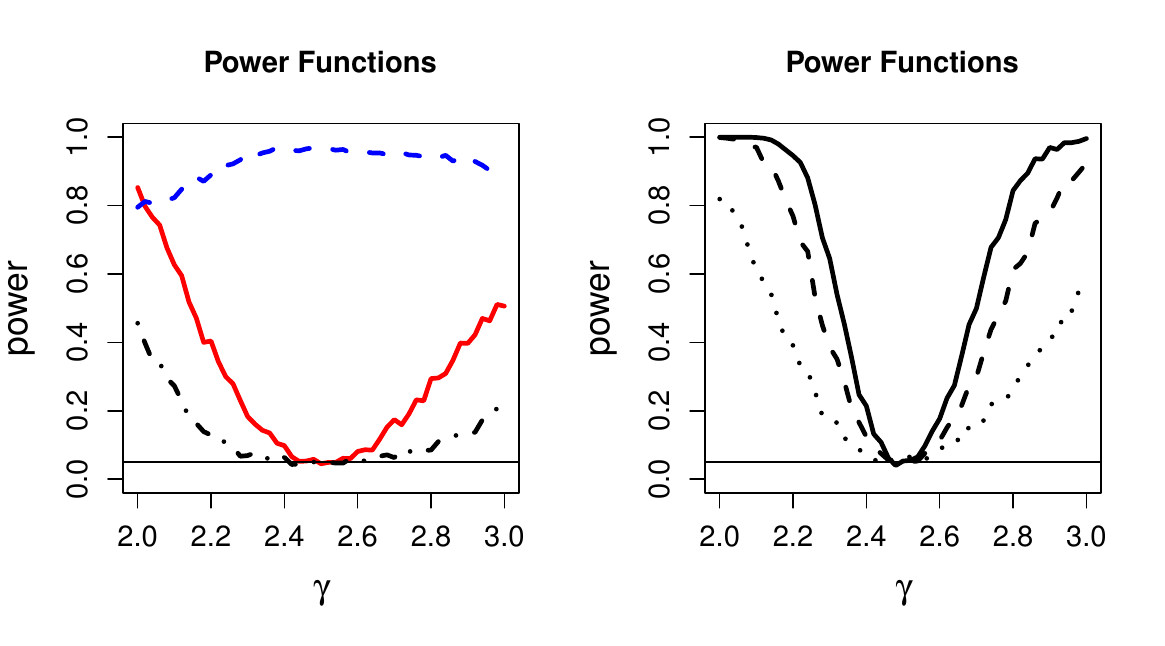}
		\caption{Empirical power functions of  $\gamma\,$ for  scale-free networks from the Barab\'{a}si-Albert model, {equipped with the Frobenius metric}, with parameter $2.5$ for $G_1$ and $\gamma\,$ for $G_2$. In the left panel, the solid red curve corresponds to the bootstrapped version of  the proposed test (\ref{rej}), the blue dashed curve to the graph based test \cp{chen:16} and the dot-dashed black curve to the energy test \cp{szek:04}. 
		In the right panel, the solid power function  corresponds to the proposed asymptotic test (\ref{rej}) for $n_1=n_2=450$, the dashed power function  to the test for $n_1=n_2=250$ and the dotted power function  to the test for $n_1=n_2=100$. The level of the tests is $\alpha= 0.05$ and is indicated by the line parallel to the $x$-axis.} 
		\label{fig: fig_4}
\end{figure}

Next we consider samples of graph Laplacians of scale free networks from the Barab\'{a}si-Albert model with the Frobenius metric. These popular networks have power law degree distributions and are commonly used for networks related to the world wide web, social networks and  brain connectivity networks. For scale free networks the fraction $P(c)$ of nodes in the network having $c$ connections to other nodes for large values of $c$ is approximately $c^{-\gamma}$, 
with $\gamma$ typically in the range $2 \leq \gamma \leq 3$. Specifically, we used the Barab\'{a}si-Albert algorithm to generate samples of scale free networks with 10 nodes, as one might encounter in brain networks.  For group $G_1$, we set $\gamma=2.5$ and for group $G_2$ we selected a fixed  $\gamma$ in the interval $2 \leq \gamma \leq 3$, studying the empirical power as a function of  
$\gamma$. The left panel in Figure \ref{fig: fig_4} indicates that in this scenario the proposed test has better power behavior  than both the graph based test and the energy test. The graph based test has a high false positive rate. The right panel in Figure \ref{fig: fig_4} shows empirical evidence  that the proposed  test is also consistent in this scenario as  sample size increases. {Additional empirical results where $n_1 \neq n_2$ and power comparisons between the asymptotic and bootstrap versions of our proposed test for smaller sample sizes can be found in the Supplementary Material.}

\begin{figure}
	\centering
	\includegraphics[width=0.85\linewidth]{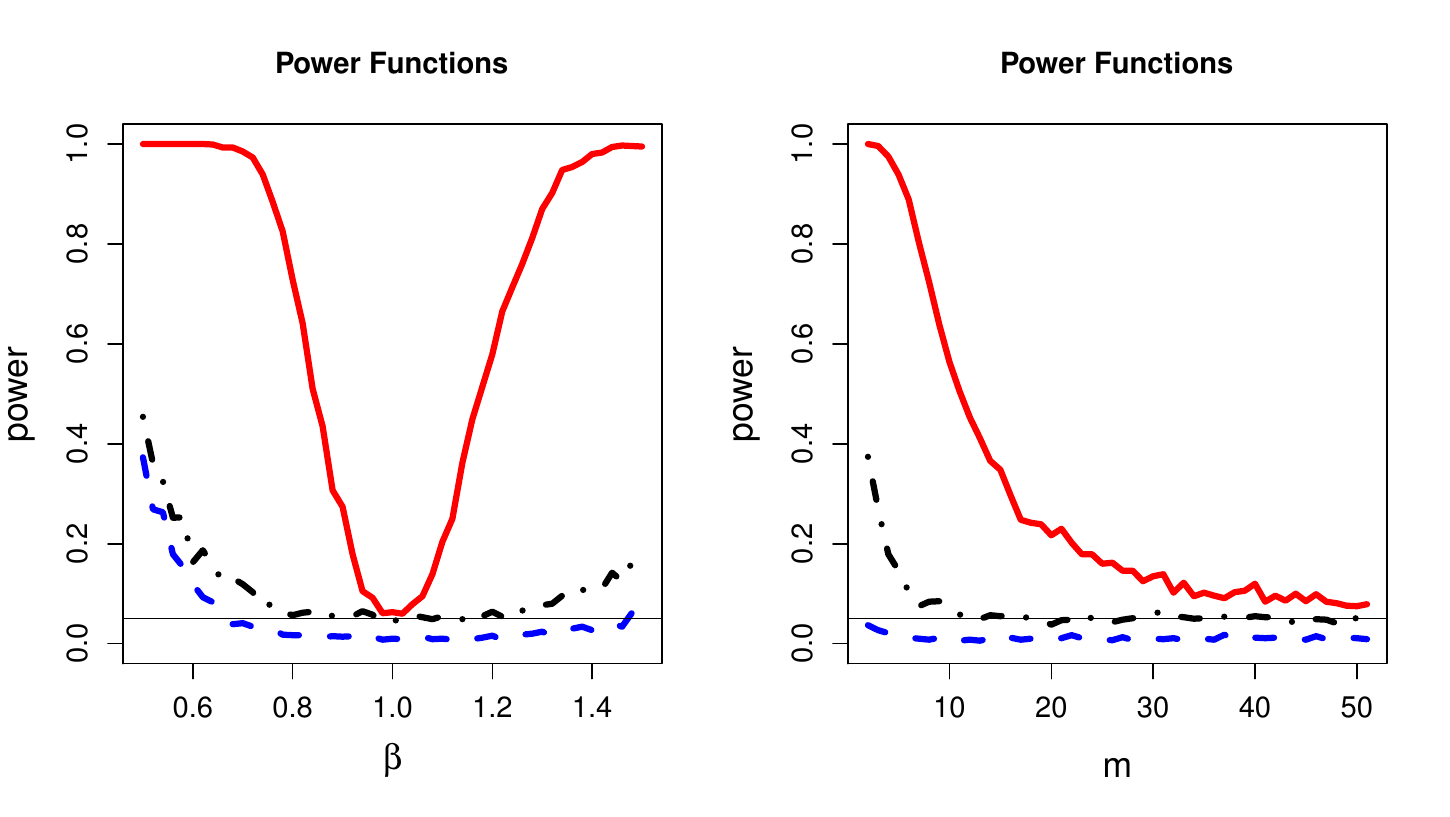}
	\caption{Empirical  power as function of  $\beta$ for  $5$-dimensional vectors, where each component  is distributed independently as $Beta(1,1)$ for group $G_1$ and as $Beta(\beta,\beta)$ with $\beta$ varying between $0.5 \leq \beta \leq 1.5 $ for group $G_2$ (left). Empirical power as function  of degrees of freedom $m$ for  $5$-dimensional vectors which are distributed independently as truncated $N(0,I_5)$ for group $G_1$ and as truncated multivariate t-distribution $t_m(0,I_5)$ with varying degrees of freedom between $2 \leq m \leq 51$ for group $G_2$ with each component of the vectors truncated to lie between $[-5,5]$ (right).  Sample sizes are $n_1=n_2=100$ and level 
		$\alpha=0.05$. 
		Solid red curves correspond to the bootstrapped version of the  proposed test(\ref{rej}), dashed blue curves to the graph based test \cp{chen:16} and dot-dashed black curves to the energy test \cp{szek:04}. The level of the test is set to 0.05 and is indicated by the line parallel to the $x$-axis. }
	\label{fig:fig_7}
\end{figure}

In the multivariate setting we considered 5-dimensional vectors. We took the five components of the vectors to be distributed independently as $Beta(1,1)$ for group $G_1$, while  for group $G_2$, the  five components were assumed to be distributed independently as $Beta(\beta,\beta)$.  Empirical power was then obtained as a function of  $\beta$ for  $0.5 \leq \beta \leq 1.5$.  In a second  scenario, the vectors are distributed as truncated multivariate normal distributions $N(0,I_5)$ for group $G_1$,   where each of the components was truncated to lie between $[-5,5]$.  For group $G_2$, we chose  a 5-dimensional t-distribution $t_{m}(0,I_5)$,  $m$ indicating  the degrees of freedom. As the degrees of freedom $m$ increase, the shape of the distribution of $G_2$   becomes more similar to that of group $G_1$. We obtained the empirical  power  as function of  $m$, for  $2 \leq m \leq 51$. Figure \ref{fig:fig_7} illustrates  that  the proposed test outperforms the comparison tests in this scenario.

\section{Data Illustrations}
\label{sec: data}
\subsection{Mortality Data}
	\no The Human Mortality Database provides data in the form of yearly lifetables differentiated by countries and gender. Presently it includes  yearly mortality data for 37 countries, available at  \url{<www.mortality.org>}. These can be converted to a density of age-at-death for each country, gender and calendar year, by first converting the available lifetables into histograms and then applying local least squares smoothing,  for which we used the Hades package at 
\url{<http://www.stat.ucdavis.edu/hades/>}  with bandwidth $h=2$.   
The random objects we consider are the resulting densities of age-at-death. Considering  the time period 1960-2009 and the  31 countries in the database for which records are available  for this time period, we obtained the densities of age at death for the age interval $[0,80]$.  

\begin{figure}[H]
	\centering
	\includegraphics[width=0.6\linewidth]{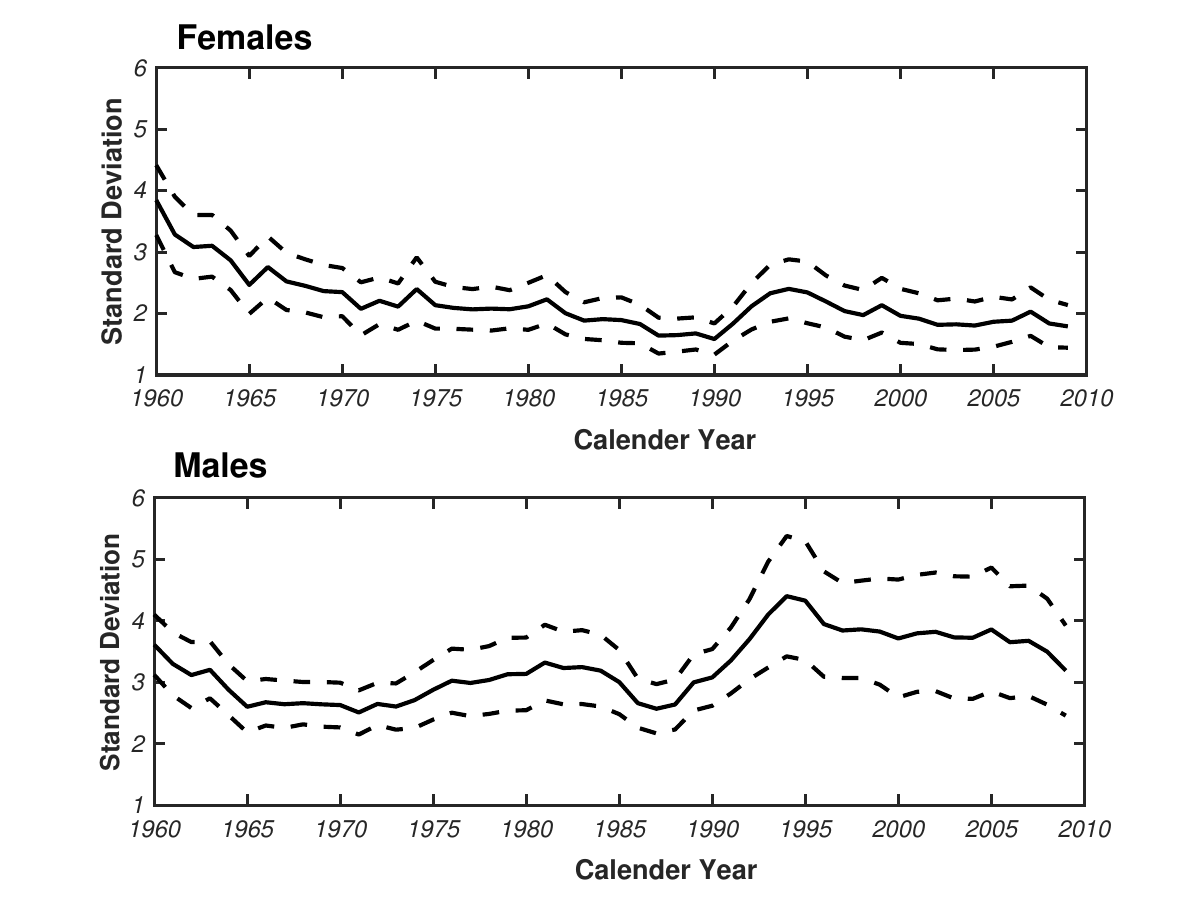}
	\caption{Yearly \F \ standard deviations (solid lines) along with 95 \% pointwise bootstrap confidence limits (dashed lines) for age-at-death densities for various countries in different calendar years, for  females (top) and males (bottom).}
	\label{fig:fig_8}
\end{figure}

From these densities we obtained 
quantile functions to compute  the Wasserstein distance, which is the metric we adopted for this distribution space.  The  \F \ standard deviations were  computed as a function of calendar year, and  are shown  along with pointwise 95 \% bootstrap confidence bands  in Figure \ref{fig:fig_8}, separately for males and females. One finds  that there is a small peak in variance of mortality between 1980-1985 for males followed by a larger peak between 1993-1996. For females, this later peak is also quite  prominent. These peaks might possibly be attributed to major political upheaval in Central and Eastern Europe during that period since several countries in the  dataset belong to these regions. 
The countries in the dataset that experienced  some turmoil associated with the end of Communist role are Belarus, Bulgaria, Czech Republic, Estonia, Hungary, Latvia, Poland, Lithuania, Russia, Slovakia and Ukraine. {The other countries in the dataset are Australia, Austria, Belgium, Canada, Denmark, Finland, France, Iceland, Ireland, Italy, Japan, Luxembourg, Netherlands, New Zealand, Norway, Spain, Sweden, Switzerland, United Kingdom and United States of America.}

\begin{figure}
	\centering
	\includegraphics[width=0.85\linewidth]{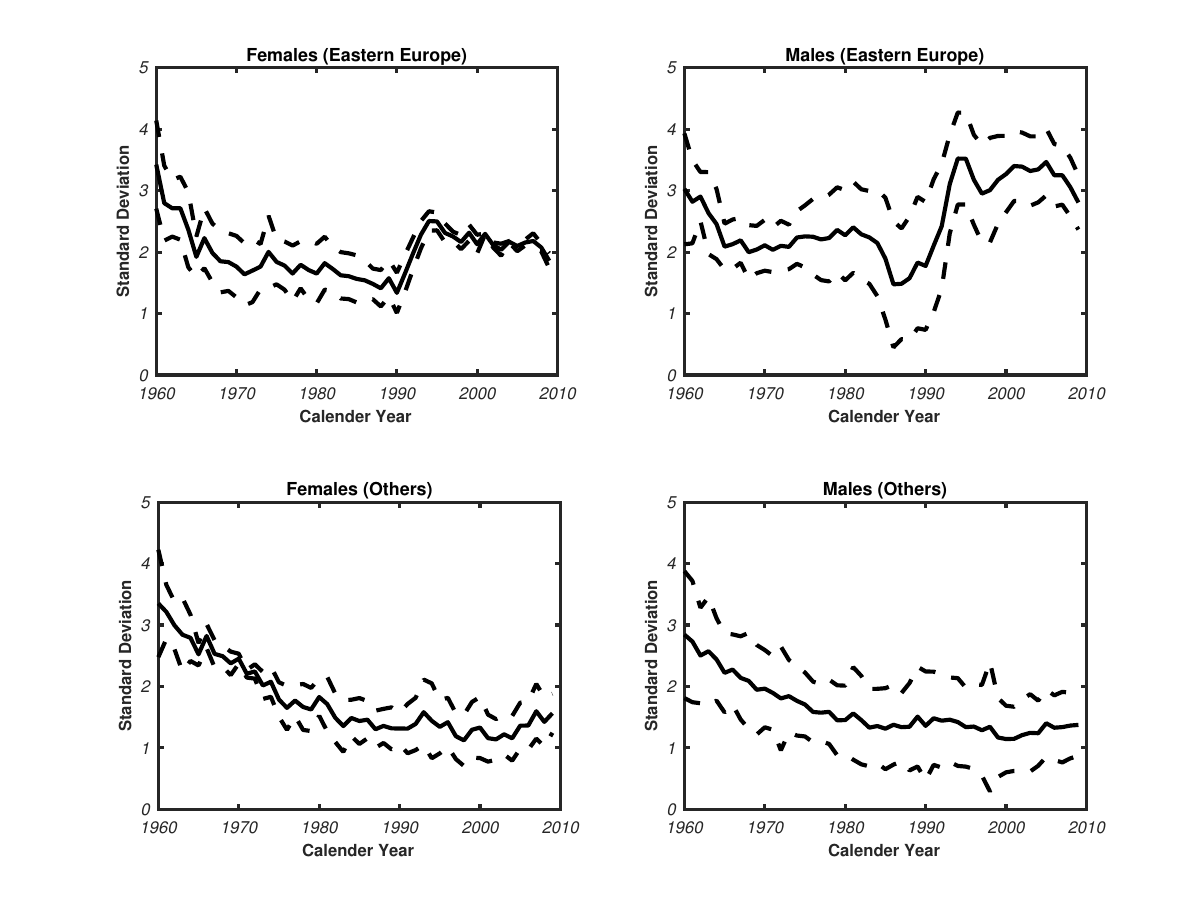}
	\caption{Yearly \F \ standard deviations (solid lines)  with 95\% pointwise bootstrap confidence limits (dashed lines) for age-at-death densities for females in Eastern European countries (top left), females  in other countries (bottom left), males in Eastern European countries (top right) and males in other countries (bottom right).}
	\label{fig:fig_9}
\end{figure}

\begin{figure}
	\centering
	\includegraphics[width=0.6\linewidth]{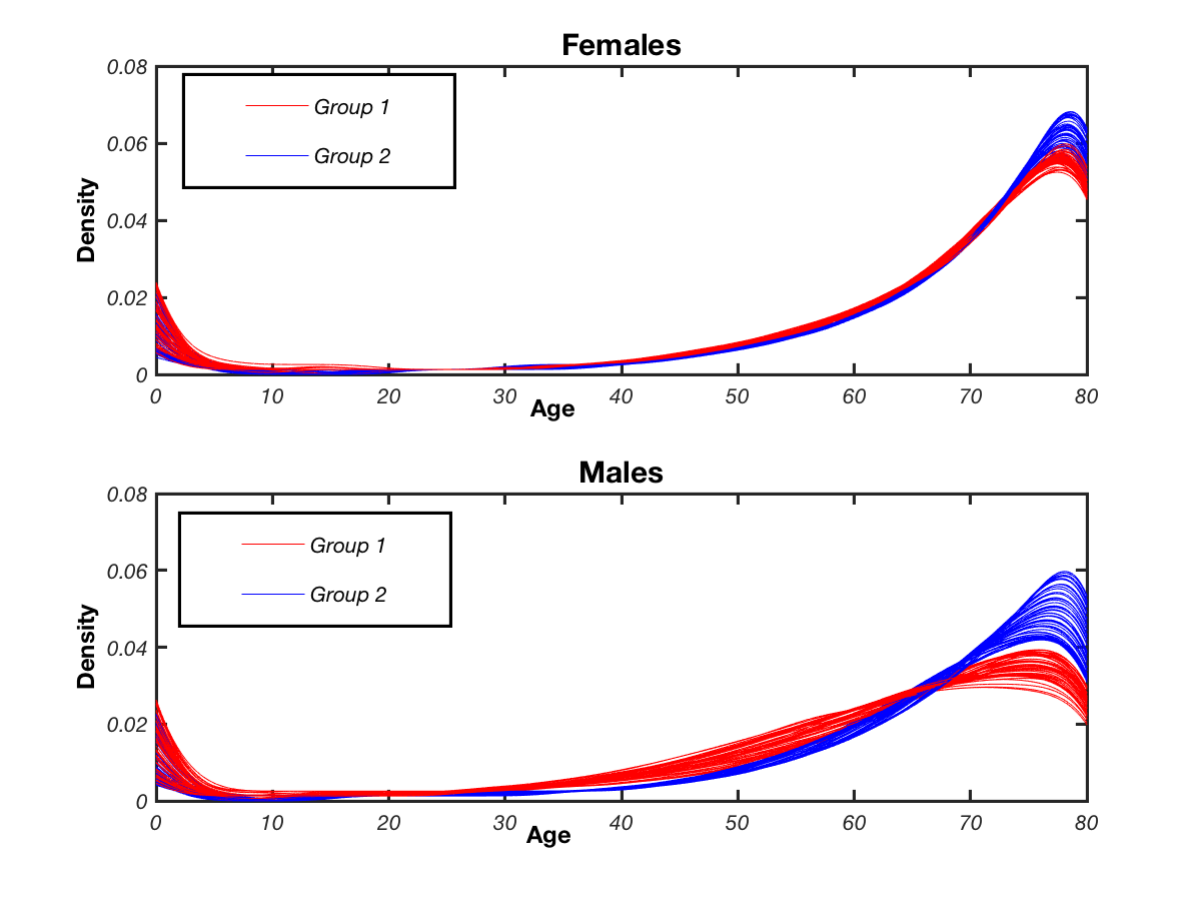}
	\caption{Wasserstein-\F \ mean age-at-death densities for  the years 1960-2009 for Eastern European countries (red) and other countries (blue), for females (top) and males (bottom). }
	\label{fig:fig_10}
\end{figure}
To visually check whether the variance peak around 1990s is indeed due to these countries, we  split our dataset into two groups, group 1 consisting  of the above Eastern European countries and group 2 of all other countries. We then created separate   plots for these two groups.  As indicated in Figure \ref{fig:fig_9},  for group 2 the variance of age-at-death distributions indeed has a decreasing trend over the years for both males and females, while the variance shows distinct fluctuations for  both males and females in group 1.  The group-wise Wasserstein \F \ mean densities of the countries for the various calendar years are illustrated in Figure \ref{fig:fig_10}. The mean densities suggest that there is a  difference between the two groups for both males and females during a significant portion of the time period between 1960 to 2009.

\begin{figure}
	\centering
	\includegraphics[width=0.6\linewidth]{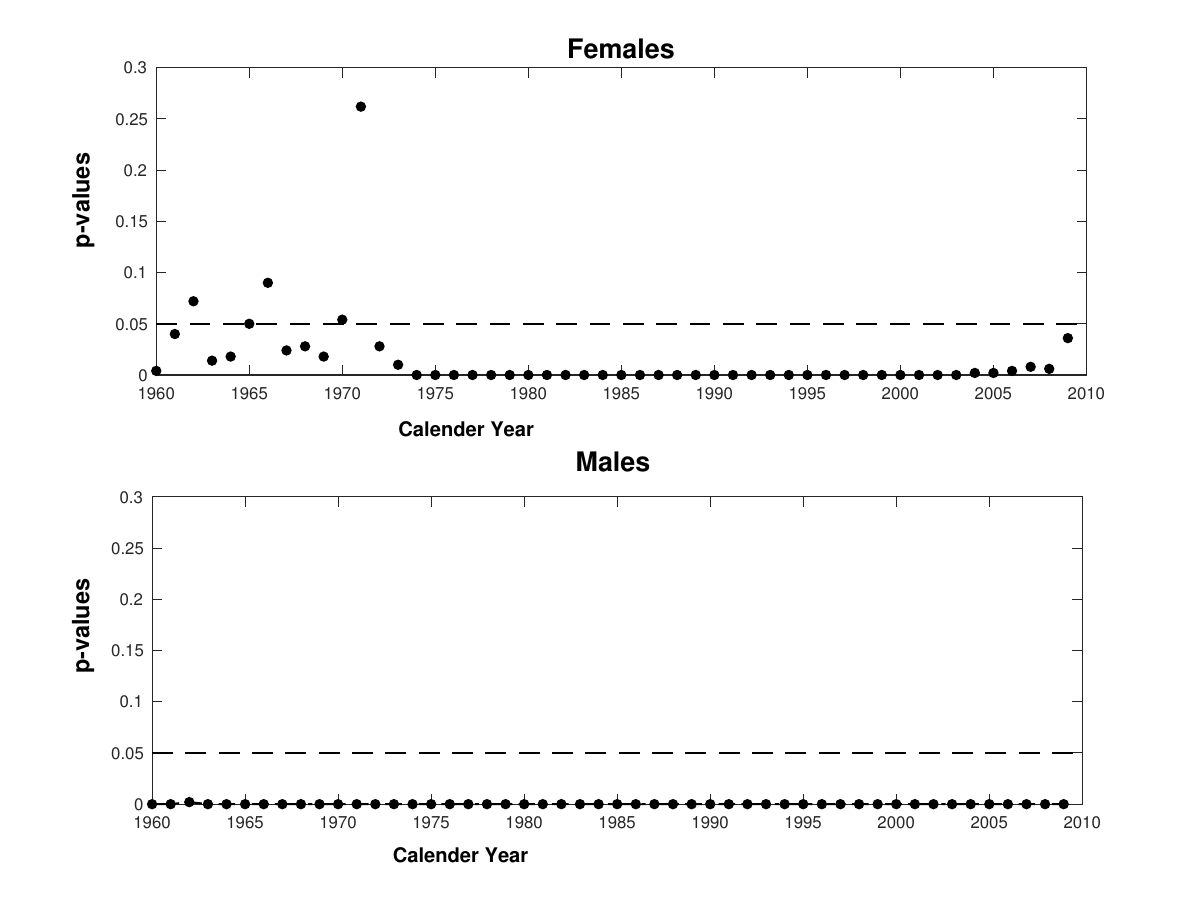}
	\caption{Dots denote the p-values for testing the differences in population age-at-death distributions between Eastern European and all other countries for each of various calendar  years for females (top) and males (bottom) with the bootstrap version of the proposed test (\ref{rej}). The dashed horizontal line indicates the .05 level.}
	\label{fig:fig_11}
\end{figure}

To more formally  test for differences in age-at-death distributions between groups 1 and 2 as defined above, we carried out the  bootstrap version of the proposed test to accommodate for the relatively small sample sizes. Figure \ref{fig:fig_11} illustrates  the $p$-values obtained for each year when applying test (\ref{rej}).  The $p$ values for testing the null hypothesis that \F\ means and variances are the same between groups 1 and 2 are far below 0.05 since 1972 for females and during the entire considered time period since 1960 for males, providing  evidence that there is a systematic difference between the Eastern European countries and the other countries in the data set during this time period.

\subsection{Comparing intra-regional connectivity for Alzheimer's and mentally normal subjects  using fMRI data}
\noindent Alzheimer’s disease  is an irreversible, progressive neuro-degenerative brain disorder that slowly destroys memory and thinking skills, eventually leading to severe dementia. Alzheimer's has been found to have associations with abnormalities in functional integration of brain regions. Recent studies as in 
\cite{sui:15} have indicated that Alzheimer's selectively targets regions of high-connectivity in the brain. Such regions are often referred to as hubs. 
The posterior midline, in particular the posterior cingulate/precuneus (PCP) brain region as described in \cite{buck:09} is a hub of high cortical connectivity and functional connectivity in this region could be a potential biomarker for Alzheimer's. For each hub region, a so-called seed voxel is identified as the voxel with the signal that has the highest correlation with the signals of nearby voxels. To quantify intra-hub connectivity, following \ci{mull:16:1},  we   analyze the distribution of the correlations between the signal at the seed voxel of the PCP hub and the signals of all other voxels within an $11 \times 11 \times 11$ cube of voxels that is centered at the seed voxel.  

After removal of data with outliers and  corrupted signals, the   subjects in our analysis consisted of cognitively normal
elderly patients and demented elderly patients diagnosed with Alzheimer's  (after removal of outliers), each of whom underwent an fMRI scan at the UC Davis Imaging Research Center. Preprocessing of the recorded BOLD (blood-oxygenation-level-dependent) signals was implemented by adopting the standard procedures of slice-timing correction, head motion correction and normalization, in addition to linear detrending to account for signal drift and band-pass filtering to include only frequencies between 0.01 and 0.08 Hz. The signals for each subject were recorded over the interval [0, 470]  (in seconds),
with 236 measurements available at 2 second intervals.

The study included 171 normal subjects but since Alzheimer's is a disease that is known to progress with age, to unambiguously rule out a differential effect of age on connectivity comparisons between the two groups,  only 87 out of the 171 subjects were included in the comparison, by matching their ages with those of the  demented patients. To check that the age matching worked, the age distributions of the included 87 normal elderly subjects  and the 65 Alzheimer's patients were compared with  the Wilcoxon rank sum test for the null hypothesis of equal age distribution across the two groups, which yielded  a $p$-value of 0.84.

\begin{figure}[H]
	\centering
	\includegraphics[width=0.5\linewidth]{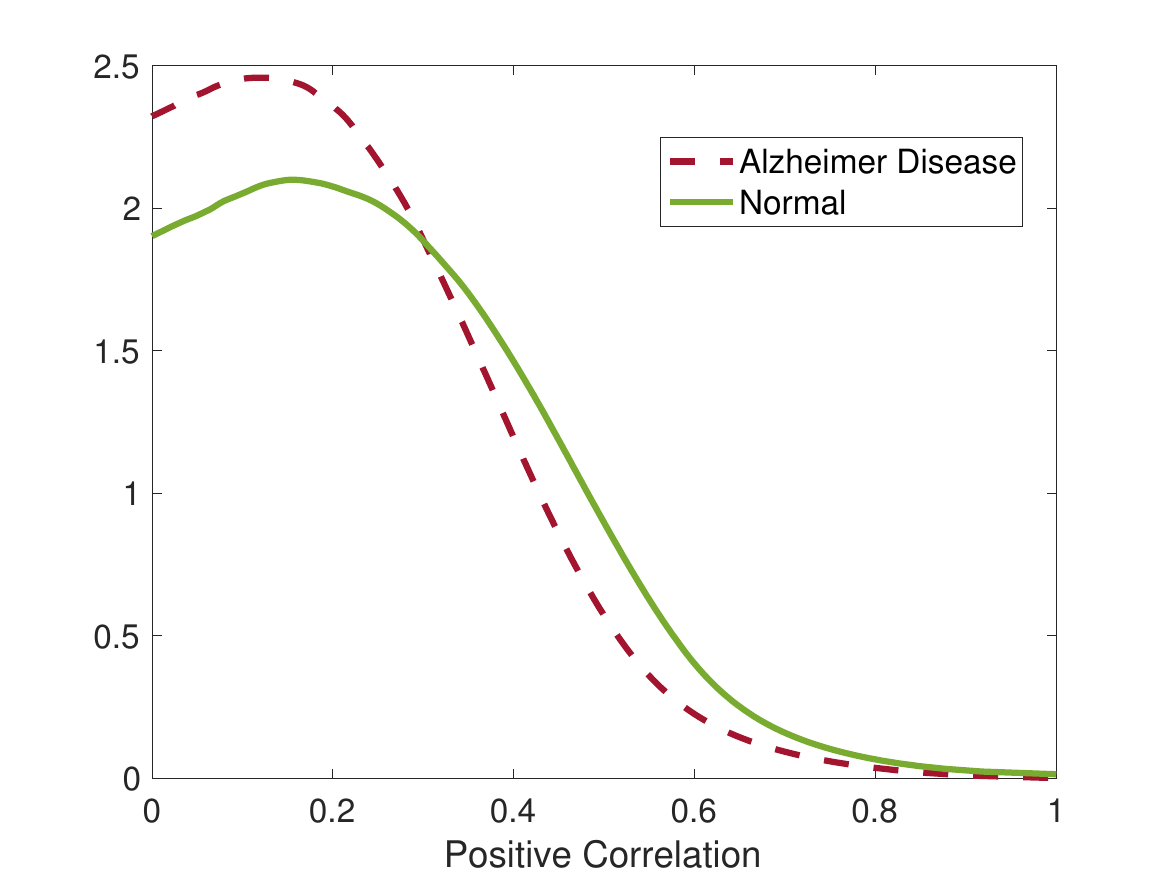}
	\caption{Wasserstein mean probability distributions of positive correlations as density functions in the PCP hub for normal subjects (green and solid) and Alzheimer disease subjects (red and dashed).}
	\label{fig:fig_13}
\end{figure}

For each subject, the target is the density function of positive correlations within the PCP hub, where this density was estimated from the observed correlations with a  kernel density estimator, utilizing the standard Gaussian kernel and bandwidth  $h = 0.08$.  As negative correlations are commonly ignored in connectivity analyses, the densities were estimated on
[0, 1]. The resulting sample of densities is then an i.i.d. sample across subjects. Figure \ref{fig:fig_13} shows the Wasserstein Fr\'{e}chet mean probability distributions represented as density functions.
To compare the two populations  of distributions, we   applied the asymptotic and the bootstrap version of the proposed test to these samples of density functions, where the asymptotic version yielded  a $p$-value of $p = 0.002$  and the bootstrap version a p-value of $p=0.001,$ indicating  that significant differences exist  in terms of intra-regional connectivity between Alzheimer's patients and age-matched normal subjects.

\subsection{Comparing brain networks of Alzheimer's patients}
\noindent It is well known  that brain hubs, being regions of high connectivity in the brain, are connected for functional integration of their specialized roles \cp{spor:11}.  Studying interconnections between hubs can reveal important insights about brain diseases like Alzheimer's. Disorders of cognition can be associated with disrupted connectivity between cortical hubs as discussed in \cite{buck:09}. One question of interest is whether the interconnections change with age in subjects with Alzheimer's.  To study this question with the proposed test, we consider connections between the 10 cortical hubs listed  in Table 3 of \cite{buck:09}, which are also referred to as regions of interest.

For this analysis, we considered the  65  subjects with Alzheimer's that were  discussed in the  preceding subsection. For each subject, a $10 \times 10$ connectivity matrix was obtained, with  entries that correspond to the observed  correlations between average fMRI signals obtained  from $3 \times 3 \times 3$ cubes around the seed voxels of the 10 hubs. The  entries of this matrix are the  so-called Pearson correlations that are commonly used in neuroimaging to quantify the correlation of random signals. These Pearson correlations can be understood as a version of dynamic correlation that has been studied previously in functional data analysis \cp{mull:05:2}.

Following standard practice in neuroimaging, these 
subject-specific connectivity matrices were then thresholded at correlation level 0.25 \cp{buck:09} to obtain adjacency matrices of networks with the hubs as the nodes, so that the presence of an edge indicates a  correlation greater than 0.25. Subject-specific graph Laplacians were then formed from these adjacency matrices.

\begin{figure}
	\centering
	\includegraphics[width=0.8\linewidth]{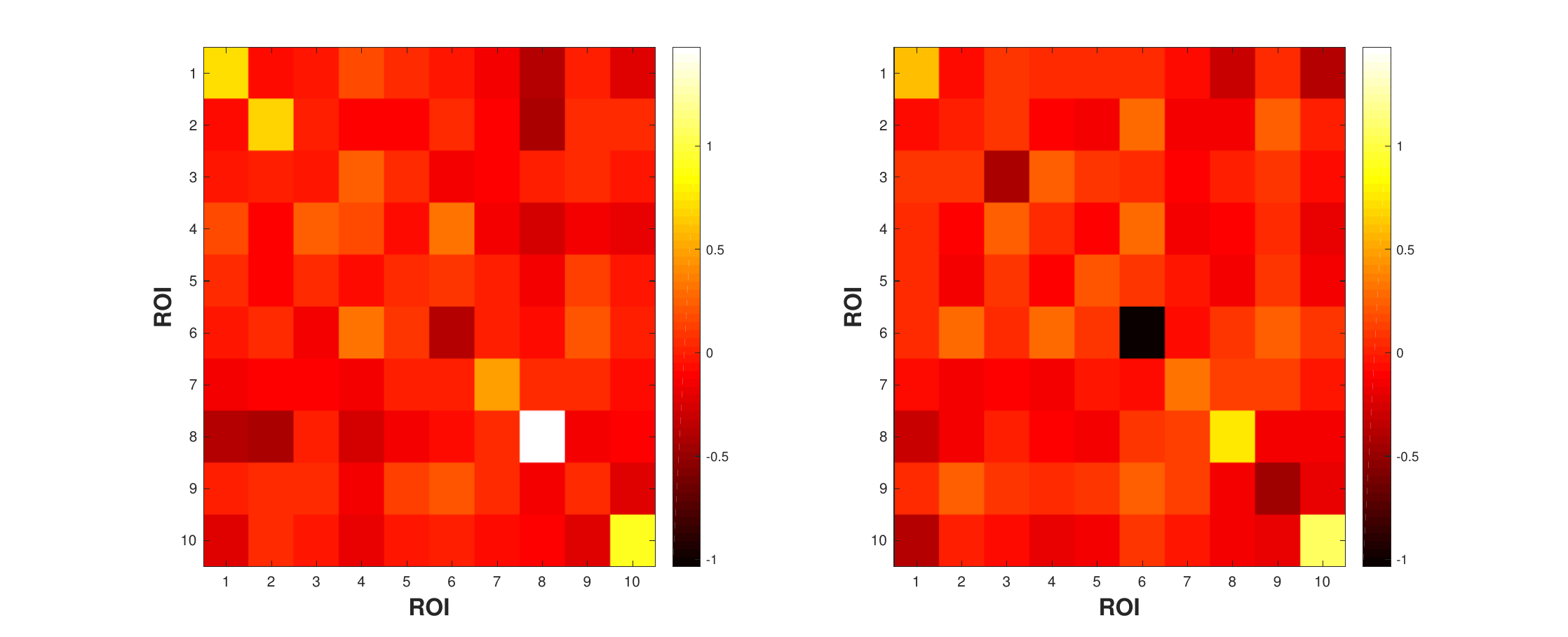}
	\caption{The left panel shows the difference between the average graph Laplacians of subjects with Alzheimer's aged between 70 and 80 and those aged 70 or below. The right panel shows the analogous comparison  for the group aged 80 or above and the group aged 70 or below. Here ROI stands for regions of interest which refer to the 10 cortical hubs that are listed  in Table 3 of \cite{buck:09}. }
	\label{fig:fig_14}
\end{figure}

These cognitively impaired patients were split into three groups based on their ages. Subjects were assigned to groups $G_1$, $G_2$ or $G_3$ based on whether they were aged $70$ or below, between age 70 and $80$ or $80$ and above. The left panel in Figure \ref{fig:fig_14} shows the difference of the average graph Laplacians of subjects in group $G_2$ and subjects in group $G_1$ and the right panel the difference of the average graph Laplacians of subjects in group $G_3$ and subjects in group $G_1$.

Since the group sample sizes of $G_1$, $G_2$ and $G_3$ are small we applied the bootstrap version of the proposed  test (\ref{rej}) to determine whether there  are  significant differences between the networks of the three age groups. The null hypothesis of equality of means and variances of the population distributions of the graph Laplacians was rejected with a bootstrap $p$-value of 0.032, providing some moderate evidence that the inter-regional brain  networks change with  age for subjects with Alzheimer's.

\section{Discussion}
We propose a  straightforward extension of analysis of variance for metric space valued data objects. Location and scale are the predominant modalities by which differences in populations of random variables in Euclidean spaces are assessed.  Differences in \F \ means and variances in the distributions of populations of random objects can be regarded as a generalization to the more general case of metric space valued objects. The proposed test is shown to have  good finite sample performance in the simulations, is free of tuning parameters and therefore can be easily applied. 
The proposed tests in both their asymptotic and  bootstrap versions are supported by asymptotic results  derived from empirical process theory and are useful in various applications.  

\section*{Supplementary Material}
Supplementary material available at {\it Biometrika} online includes the proofs of main and auxiliary results, justification for the bootstrap version of the test and additional simulations. 

\section*{Acknowledgements}
We wish to thank the reviewers for constructive and most helpful comments that led to various corrections and improvements. This research was supported by the National Science Foundation.

\references
\ed